\documentclass[12pt]{article}

\usepackage{amsmath,amssymb,latexsym}
\usepackage{booktabs}
\usepackage{color}
\usepackage[round]{natbib}

\newtheorem{theorem}{Theorem}[section]

\newtheorem{lemma}[theorem]{Lemma}

\newtheorem{corollary}[theorem]{Corollary} 
\newtheorem{remark}[theorem]{Remark}

\newtheorem{definition}[theorem]{Definition}

\def\lcm{{\rm lcm}}

\newcommand{\GL}{\mathrm{GL}}
\newcommand{\SL}{\mathrm{SL}}
\newcommand{\SU}{\mathrm{SU}}
\newcommand{\GU}{\mathrm{GU}}
\newcommand{\Sp}{\mathrm{Sp}}

\newcommand{\SO}{\mathrm{SO}}

\newcommand{\ppd}{\mathrm{ppd}}
\newcommand{\full}{\mathrm{full}}
\newcommand{\F}{\mathbb{F}}

\title
{Elements in finite classical groups whose powers have large 1-Eigenspaces}

\author{Alice  C. Niemeyer$^{1,2}$, Cheryl E. Praeger$^{1,3}$}
\begin{document}

\maketitle

\begin{abstract}
We estimate  the proportion of  several classes of elements  in finite
classical groups which are readily recognised algorithmically, and for
which some power has a large fixed point subspace and acts irreducibly
on a complement of it.   The estimates are used in complexity analyses
of new recognition algorithms for finite classical groups in arbitrary
characteristic.
\end{abstract}


\section{Introduction}
A crucial task  in designing algorithms to compute  with matrix groups
defined  over finite  fields is  to recognise  whether a  given matrix
group $H$  is isomorphic to  a finite classical  group and, if  so, to
find  an  isomorphism with  the  natural  representation  of $H$.  The
isomorphism is defined by identifying appropriate elements of $H$ with
the  standard  generators of  the  natural representation.  Algorithms
which accomplish  this task are generally referred  to as constructive
recognition algorithms by \cite{KS}. In this paper, we concentrate on the
important  special case  when  $H$  is already  given  in its  natural
representation, by a set of $d  \times d$ matrices over a finite field
$\F_q$ (but we have not yet found appropriate elements of $H$ that can
serve as  standard generators). In this  special case, for  odd $q$, a
highly  efficient  and  practical   algorithm  has  been  designed  by
\cite{LGO}.  Recently, \cite{Letal} 
designed recognition algorithms for even  $q$ and \cite{NS2} (see also
\cite{NS}) designed algorithms for arbitrary $q$. 
 In both new
algorithms,  a pivotal  step  is to  find  an element  $x$ which  acts
irreducibly on a subspace of dimension $m$, that is it does not leave
invariant a subspace of the $m$-dimensional space,  with $m$ small relative to
$d$, and fixes a complementary subspace pointwise. 
\cite{Letal} seek such elements elements for which the dimension $m$
is at least a constant fraction of $d$ 
(see Section~5, especially Lemma~5.5 of their paper),
whereas \cite{NS2} use such elements with $m = O(\log(d))$ (see our
Corollary~\ref{cor:main3}).
Using results of \cite{PSY},  \cite{NS2} prove 
that, if $x$ satisfies an  additional condition on its order (divisible
by a primitive prime divisor, see Section~\ref{sec:Q}) then, with high
probability,  $x$ and  a random  conjugate of  $x$ in  $H$  generate a
classical group  of dimension $2m$. This  smaller classical group
is then processed in a recursive procedure.
As elements  with large $1$-eigenspaces are rare,  the algorithms seek
such elements  as appropriate  powers of nearly  uniformly distributed
random elements  of $H$.   In this paper,  we show that  a significant
proportion of elements of $H$ have a power with the desired properties.

For  an  integer  $k$  we introduce  in  Definition~\ref{def:Q}  three
subsets of $H$, called  $Q_k, Q_k^{\ppd}$ and $Q_k^{\full}$, for which
$Q_k^\full\subseteq  Q_k^\ppd\subseteq  Q_k$  and  if  $k\geq3$,  then
$Q_k^\full\ne\emptyset$.  
(Here $H$ is one of the classical groups in Table 1, and the defining
condition for the sets depends on the integer $\alpha$ given in
Table~\ref{tab:G}.)
Elements in $Q_k$ power to elements
leaving invariant an $\alpha k$-dimensional subspace without a 1-eigenvalue
and having a $(d-\alpha k)$-dimensional $1$-eigenspace
and in addition elements in $Q_k^\ppd$ 
power to elements acting irreducibly on the $\alpha k$-dimensional subspace.
Membership in these sets can readily be  determined algorithmically: for example, in
the case of $Q_k$, by examining the degrees of the irreducible factors
of  the characteristic  polynomial.
Properties  of these  subsets are
discussed  in  more   detail  in  Section~\ref{sec:Q}.  
We  state  here   our  main  result  on  the
properties and proportions of elements in these subsets.  We note\
that  the proportion of  elements of  $H$ that  lie in  $Q_k^\ppd$ for
$\alpha k  > d/2$ differs  only by a  constant from the  proportion of
so-called $\ppd$-elements  determined in \cite[Theorem~5.7]{NP98} (see
also its generalisation in \cite[Theorem~1.9]{NP}). Throughout the
paper $\log$ denotes the natural logarithm to base $e$.

\begin{theorem}\label{the:short-main}
Let $H$ be a $d$-dimensional classical group  defined over a finite 
field with $q$ elements, as in one of the rows of Table~$\ref{tab:G}$.
Let $n,\alpha$ be as in Table~$\ref{tab:G}$, and let $k$ be an  
integer such that  $\max\{3,\log(n)\}\le k\le n/2$ and $k$ odd if $H=SU_n(q)$. 
Let $Q$ denote one of the sets $Q_k$, $Q_k^{\ppd}$ or $Q_k^{\full}.$  
Then each element of $Q$ powers up to an 
element which has a $(d-\alpha k)$-dimensional $1$-eigenspace and 
acts irreducibly 
on a complementary invariant subspace of
dimension $\alpha k$. Moreover for $Q = Q_k$ or $Q=Q_k^\ppd$
\begin{displaymath}
\frac{2}
{9 e\alpha k} \le \frac{|Q|}{|H|} \le 
\frac{5}{3\alpha k} \mbox{\quad while\quad }
\frac{2}
{10 e\alpha k^2\log(q+1)} \le \frac{|Q_k^\full|}{|H|} \le 
\frac{5}{3\alpha k}.
\end{displaymath}
\end{theorem}

\begin{table}[t]
\begin{center}
\begin{tabular}{lllclc}
\toprule
$H$ & $d$  & $\alpha$& $\delta$ & $\epsilon$\\ 
\midrule
$\SL_n(q)$ & $n$ &  $1$ & $1$ & $1$\\ 
$\SU_n(q)$ & $n$ & $1$  & $2$ & $-1$ \\ 
$\Sp_{2n}(q)$ & $2n$ & $2$& $1$ & $1$\\ 
$\SO_{2n+1}(q)$ & $2n+1$ & $2$& $1$ & $1$\\ 
$\SO^\pm_{2n}(q)$ & $2n$ &$2$&  $1$ & $1$\\ 
\bottomrule
\end{tabular}
\end{center}
\caption{Groups and constants in Theorems~1.1 and 3.3}\label{tab:G} 
\end{table}

Every  element $g$ in  a finite  classical group  admits a 
multiplicative Jordan decomposition as $g=us,$ where $s$ is semisimple
and  $u$ is  unipotent. As  $s$ is  diagonalisable over  the algebraic
closure $\overline{\F}_q$  of the field $\F_q$ with  $q$ elements, the
eigenvalues of  $g$ over  $\overline{\F}_q$ are those  of $s$  and the
characteristic   polynomial  of   $g$   equals  that   of  $s$.    The
characteristic  polynomial   of  $s$  yields   information  about  the
eigenvalues of  $s$ (and $g$):  for each irreducible factor  of degree
$k$,  there are $k$  eigenvalues over  $\overline{\F}_q$ which  lie in
$\F_{q^k}$ and no smaller field.

Our aim is  to group the eigenvalues of $g$ into  two disjoint sets in
such a way that we can find  an integer $B$ such that the $B$-th power
of each eigenvalue in the first  set is equal to $1$, while the $B$-th
power of  each eigenvalue in  the second set  is distinct from  1.  We
identify elements  $g$ in classical groups  for which we  can read off
from  the characteristic  polynomial  whether or  not  they have  such
eigenvalues. Moreover the integer $B$ can then also be determined from
the characteristic polynomial, see Remark~\ref{rem:B}.

In order to give a precise description of  the elements we seek, we  use
detailed information about the eigenvalues of elements in finite
classical groups. Such information was obtained in
\cite[Sections~3,5]{LNP} and we present a summary of the information
relevant to this 
paper in Section~\ref{sec:ew}. In Section~\ref{sec:Q} we are then able
to define the elements we seek precisely. 
We state in Theorem~\ref{the:main} a complete and detailed
version of Theorem~\ref{the:short-main} giving upper and lower bounds
on the proportions for the three classes of elements separately. In  Section~\ref{sec:ingr}
we outline the strategy of our proof and 
prove Theorem~\ref{the:main} in Subsection~\ref{sub:main}.
 Theorem~\ref{the:short-main} then follows.

\section{Eigenvalues of elements in maximal tori}\label{sec:ew}

The theory of algebraic groups  necessary to understand the results in
this paper can be found in \cite{Carter85}.  Here we just mention that
for a  finite classical group $H$  and for $W$  the corresponding Weyl
group,  the $H$-conjugacy classes  of $F$-stable  maximal tori  are in
one-to-one correspondence  with the $F$-conjugacy classes  of the Weyl
group $W$ of $H$.  When $H$ is a finite classical group the Weyl group
$W$ of $H$  is isomorphic to a full symmetric group  $S_n$ (lines 1-2
of Table~\ref{tab:G}), or a wreath
product $S_2\wr S_n$ (lines 3-4 of Table~\ref{tab:G}), or the index 2 subgroup of $S_2\wr S_n$ consisting 
of even permutations (line~5 of Table~\ref{tab:G}). Thus the
$F$-conjugacy classes of the Weyl group 
are parametrised by  partitions or  pairs of partitions of
$n$.   In this section we describe  the eigenvalues over
$\overline{\F}_q$   and   the   characteristic  polynomials   of   the
semi-simple elements  in the maximal  tori $T_C$ corresponding  to the
$F$-conjugacy  class $C$  of  $W$.  The  information  needed for  this
description  is extracted  from \cite[Sections~3,5]{LNP}.   
Let $H, n$ and $\varepsilon$ be as in one of the lines of
Table~$\ref{tab:G}$.

\subsection{Characteristic polynomials in corresponding tori in 
$\SL_n(q)$ and $\SU_n(q)$}

Here $H$ is the set of elements of the algebraic group $G = \GL_n(\overline{\F}_q)$ 
fixed by the Frobenius morphism $F:(a_{ij}) \rightarrow (a_{ij}^{\epsilon q})$, 
where $\epsilon = 1$ or $-1$, according as $H=\SL_n(q)$ or $\SU_n(q)$,
respectively. An $F$-conjugacy class $C$ of $W$ is
uniquely determined by  a partition of $n$. For each part $\lambda$ of
this  partition the elements  of $T_C$  have $\lambda$ corresponding eigenvalues in
$\overline{\F}_q$ which  lie in $\F_{q^\lambda}$ and  these are permuted
by the map $\tau: a\mapsto a^{\epsilon q}.$ If one of
these   $\lambda$   eigenvalues  lies   in   a   smaller  field,   say
$\F_{q^{\lambda'}}$, for $\lambda'$ a proper divisor of $\lambda$, then
so do all of its Galois conjugates. Therefore, for each part $\lambda$
of the partition  of $n$ and for all  divisors $\lambda'$ of $\lambda$,
the  torus $T_C$  contains elements  whose  characteristic polynomials
have  $\lambda/\lambda'$  factors   of  degree  $\lambda'$  which  are
irreducible over $\F_q$.

If $\lambda'$ is odd, then each of these irreducible factors of degree
$\lambda'$ is also irreducible  over $\F_{q^2}$, whereas if $\lambda'$
is even then each irreducible factor of degree $\lambda'$ is a product
of a pair of $\tau$-conjugate factors of degree $\lambda'/2$ which are
irreducible over $\F_{q^2}$.

\subsection{Characteristic polynomials in corresponding tori in $\Sp_{2n}(q)$ and
  $\SO_{2n+1}(q)$} 

Here  $H$  is  the  set  of  elements of  the  algebraic  group  $G  =
\Sp_m(\overline{\F}_q)$  or   $\SO_m(\overline{\F}_q)$  (subgroups  of
$\GL_m(\overline{\F}_q)$) fixed by  the Frobenius morphism $F:(a_{ij})
\rightarrow (a_{ij}^{q})$, where $m = 2n$ if $H=\Sp_m(q)$ and $m=2n+1$
if $H=\SO_m(q)$.  An $F$-conjugacy class $C$ of $W$ is determined by a
partition  of  $n$, where  each  part  is  marked either  positive  or
negative.

For any positive part $\lambda,$  the elements of $T_C$ have $\lambda$
eigenvalues  in  $\overline{\F}_q$,  and  they,  together  with  their
$\lambda$ inverses,  lie in $\F_{q^\lambda}$  and are permuted  by the
map $a\mapsto  a^{q}$ in  two cycles of  length $\lambda$.  If  one of
these   $2\lambda$  eigenvalues   lies   in  a   smaller  field,   say
$\F_{q^{\lambda'}}$,  for $\lambda'$  a proper  divisor  of $\lambda$,
then so  do all of  its Galois conjugates  and all of  their inverses.
Therefore, for  each positive part  $\lambda$ of the partition  of $n$
and for all divisors $\lambda'$  of $\lambda$ the torus $T_C$ contains
elements  whose  characteristic  polynomials  have  $\lambda/\lambda'$
corresponding  pairs  of  factors,   each  of  degree  $\lambda'$  and
irreducible over  $\F_q$ and if $\xi\in\bar\F_q$  is a root  of one of
these factors, then $\xi^{-1}$ is a root of its paired factor.

For  any  negative  part  $\mu,$  the elements  of  $T_C$  have  $\mu$
eigenvalues in  $\overline{\F}_q$, and they together  with their $\mu$
inverses lie in $\F_{q^{2\mu}}$, and are permuted by the map $a\mapsto
a^{q}$ in a cycle of length $2\mu$.  Therefore, for each negative part
$\mu$ of the partition of $n$ and for all divisors $\mu'$ of $\mu$ the
torus  $T_C$ contains elements  whose characteristic  polynomials have
$\mu/\mu'$   corresponding  factors  of   degree  $2\mu'$   which  are
irreducible over  $\F_q$.  Moreover, if $\xi\in\bar\F_q$ is  a root of
one of these factors, then so is $\xi^{-1}.$

\subsection{Characteristic polynomials in corresponding tori for $\SO_{2n}^\pm(q)$}

Finally,  we  treat  this  case  similarly to  the  previous  case  by
considering    $\SO^\pm_{2n}(q)$   as    a    natural   subgroup    of
$\SO_{2n+1}(q)$. Then the maximal  tori of $\SO^\pm_{2n}(q)$ are those
tori of $\SO_{2n+1}(q)$ that  correspond to pairs $(\beta^+, \beta^-)$
of  partitions  of  $n$  with  $|\beta^+|  +  |\beta^-|=n$  such  that
$\beta^-$  has an even  number of  parts for  $\SO^+_{2n}(q)$, whereas
$\beta^-$ has  an odd number  of parts for $\SO^-_{2n}(q)$.   The Weyl
group   $W$   has    index   2   in   the   Weyl    group   $W_B$   of
$\SO_{2n+1}(\overline{\F}_q)$. For each $F$-conjugacy class $C$ in $W$
which  does  not  correspond to  a  single  part  of length  $n$,  the
characteristic  polynomials that  arise are  the same  as in  the case
$\SO_{2n+1}(q)$  (and all  the classes  $C$  we consider  are of  this
type).

\section{The elements we seek}\label{sec:Q}

Having obtained a description of the eigenvalues and characteristic
polynomials of semisimple elements, we are now in a position to
define the elements we seek. Let $H, d, n, \alpha, \delta, 
\varepsilon$ be as in one of the rows of Table~\ref{tab:G}. 
As mentioned in the introduction, we aim
to find elements whose eigenvalues belong to one of two distinct sets.
For a fixed integer $k$, we would like to find elements with $\alpha k$ eigenvalues
that lie in $\F_{q^{\delta\alpha k}}$ and no proper subfield, such that the
remaining eigenvalues all have order dividing an integer $B$ which
is divisible by no primitive prime divisor of $q^{\delta\alpha k}-1$.
(A primitive prime divisor of $q^m-1$ is a
prime $r$ dividing $q^m-1$ but not dividing $q^i-1$ for integers $i$ with $1
\le i < m.$)

We say that two polynomials $f_1(x)$ and $f_2(x)$ defined over a
finite field $\F_q$ with $q$ elements are \emph{conjugate}
if for any  root $a$ of $f_1(x)$ over the algebraic closure $\overline{\F}_q$
of $\F_q$, the element $a^{-1}$ 
is a root of $f_2(x)$. A polynomial $f(x)$ is \emph{self-conjugate} if
$f(x)$ is conjugate to itself. For an integer $k$, we can now define the set $Q_k$ and its two subsets $Q_k^\ppd$ and $Q_k^\full$.

\begin{definition}\label{def:Q}

Let $k$ be an integer such that $1\leq k\leq d$. 
\begin{enumerate}
\item[$(i)$] Let
$Q_k$ denote  the set  of all elements  $g \in H$  with characteristic
polynomial $c_g(x)$ for which one of the following conditions holds:

\begin{enumerate}
\item[$(a)$] $H = \SL_n(q)$ or $\SU_n(q)$ and   $c_g(x)$ has
an irreducible factor
of degree $k$ over $\F_q$, where $k$ is odd  in case $\GU_n(q)$, and no other
irreducible factors of degree divisible by  $k$,
\item[$(b)$] $H = \Sp_{2n}(q)$,
$\SO^\pm_{2n}(q)$ or $\SO_{2n+1}(q)$ and $c_g(x)$ has
a self-conjugate irreducible factor of degree $2k$,  no other
self-conjugate irreducible factors of degree divisible by $2k$ and no 
non-self-conjugate irreducible factors of degree divisible by $k$.
\end{enumerate}

\item[$(ii)$] Let $Q_k^{\ppd}$ denote the subset of $Q_k$ consisting of those 
elements $g$ for which, in addition, there is  a primitive prime divisor $r$ of
$q^{\delta \alpha k}-1$ which divides $|g|$. 

\item[$(iii)$] Let $Q_k^{\full}$ denote the 
subset of $Q_k^{\ppd}$ of elements $g$ such that,  for every primitive
prime  divisor  $r$  of  $q^{\delta\alpha  k}-1$ which  divides  $|H|$, the order
$|g|$ is divisible  by the $r$-part of $q^{\delta\alpha k}-1$ (that is, by the highest 
power of $r$ dividing $q^{\delta\alpha  k}-1$).

\end{enumerate}
\end{definition}

Then for an element $g\in Q_k$  we can find an integer $B$ as
described in the following remark such that $g^B$ has a
$(d-\alpha k)$-dimensional $1$-eigenspace and a $k$-dimensional  invariant
subspace on which $g$ acts irreducibly.

\begin{remark}\label{rem:B}{\rm
Let $g\in Q_k$.
We define  a positive integer
$B$ as follows
(see also \cite[Section~2.2]{LGO}).
Let $c_g(x)$ denote the characteristic polynomial of $g$ and suppose that $c_g(x)$
factors into irreducibles as $c_g(x) = \prod_{i=1}^s f_i(x)^{n_i}$,
where the degree of $f_i(x)$ is $k_i$, $k_1=\alpha k$, and $n_1=1$
(and $f_1$ is self-conjugate in the symplectic and orthogonal cases).
Let $B = p^\beta \prod_{i=2}^s  (q^{\delta k_i}-1)$, where $p$ is the
characteristic of $\F_q$ and $\beta = \lceil\log_p(\max(n_i)) \rceil.$
All eigenvalues of $g$ over $\overline{\F}_q$ corresponding to an
irreducible factor $f_i(x)$ of degree $k_i$, for $i > 1$, lie in 
$\F_{q^{k_i}}$ and hence are powered to $1$ by $B$. For $i>1$ we have  
$\gcd(q^{\delta\alpha k}-1,  q^{\delta k_i}-1) = q^{\delta\gcd(\alpha k, k_i)}-1$ 
and since 
$\gcd(\alpha k,k_i) \not= \alpha k$, 
it follows that 
$\gcd(q^{\alpha k}-1, B)$ divides $\lcm\{ q^\ell-1\,|\,\ell\ \mbox{a 
proper divisor of $\alpha k$}\}$. 
Therefore the eigenvalues of $g$ lying in $\F_{q^{\alpha k}}$, but in no 
proper subfield, are not powered to 1 by $B$. 
Thus  $g^B$ has a  $(d-\alpha k)$-dimensional 1-eigenspace and leaves
invariant a subspace $U$ of the underlying vector space of dimension
$\alpha k$ on which $g^B$ acts without a 1-eigenvector. 
Note that $g^B|_U$ may not be irreducible: for example, if $q=\alpha=2, k=3$, then 
there is an $f_1$ such that $g^B|_U$ has order $9$; under our assumptions 
it is possible to have $k_2=4$ and hence for $B$ to be divisible by 3 so that
$g^B$ has order 3 and acts reducibly on the $6$-dimensional subspace $U$.

Even though $g$ is irreducible on $U$, the element $g^B$ may not be. 
For example, if  $H=\SL_{11}(2)$, there exists $g\in Q_6$ with
$c_g(t)$ a product of irreducibles of degrees 2,3 and 6 and inducing a 
Singer cycle on a 6-dimensional space $U$, while $B=21$ and $g^B$
leaves invariant several 2-dimensional subspaces of $U$.  
For some applications  we require $g^B$ to be  irreducible on $U$, and
this  is  the  reason   for  defining  the  subsets  $Q_k^{\ppd}$  and
$Q_k^{\full}$  of $Q_k$. Indeed,  if $g\in  Q_k^{\ppd}$ then  $|g|$ is
divisible  by some primitive  prime divisor  $r$ of  $q^{\delta \alpha
  k}-1$, and $r$ does not  divide $q^{\delta k_i}-1$ for $i=2, \ldots,
s$ as  $k_i \not\equiv 0\pmod{\alpha  k}$. Thus $r$ divides  the order
$|g^B|$  and   hence  $g|_U$  is  irreducible.    Moreover,  if  $g\in
Q_k^\full$,  and $r^a$ divides  $q^{\delta \alpha  k} -1$,  then $r^a$
also divides $|g^B|$, for the same reason.  }
\end{remark}

Note that $Q_k^\full\ne\emptyset$  provided $q^{\delta\alpha k}-1$ has
a primitive  prime divisor dividing  $|H|$, and this is  almost always
true if $3\leq  k\leq n/2$ (with a small number  of exceptions such as
$k=6, H=\SL_{11}(2)$, see \cite{Zsigmondy}),  hence our need to assume
the subsets are non-empty in Theorem~\ref{the:main}.  Our main theorem
is as follows. Recall that we use natural logarithms.
\begin{table}[t]
\begin{center}
\begin{tabular}{l|ccc}
\toprule
$Q$      & $Q_k$ & $Q_k^{\ppd}$ & $Q_k^{\full}$\\
\midrule
$\ell_{k,Q}$ & $1-\frac{2}{q^{k/2}}$ &  $1-\frac{1}{\alpha k}$ & $\frac{\log(2) }{k\log(q+1)}$\\
$m_{n,Q}$    & $1-\frac{2}{q^{\log(n)/2}}$ & 
$1-\frac{1}{\alpha \log(n)}$ & $\frac{\log(2)}{2\log(n)\log(q+1)}$\\ 
\bottomrule
\end{tabular}
\end{center}
\caption{Definitions of $\ell_{k,Q}$ and $m_{n,Q}$ in Theorem~\ref{the:main}}\label{tab:ell} 
\end{table}

\begin{theorem}\label{the:main}
Let $q$ be a prime power and $\F_q$ a finite field with $q$ elements. 
Let $H$, $d$, $\alpha, \delta$ and $\epsilon$ be as in one of the lines
of Table~$\ref{tab:G}$. Let 
$k$ be a positive integer with $\max\{3,\log(n)\} \le k\leq n/2$ 
and such that   
$k$ is odd when $H = \SU_n(q)$. 
Let $Q$ be one of $Q_k, Q_k^{\ppd}, Q_k^{\full}$ such that $Q\ne\emptyset$,
and let $\ell_{k,Q}$ and $m_{n,Q}$ be as in Table~$\ref{tab:ell}$. Then 
\begin{enumerate}
\item[$(a)$]
\begin{displaymath}
\frac{\ell_{k,Q}}{3e\alpha k} \le \frac{|Q|}{|H|} \le 
\frac{5}{3\alpha k}.
\end{displaymath}
\item[$(b)$]
If, in addition, 
$\log(n) < k \le 2 \log(n)$, then
\begin{displaymath}
\frac{m_{n,Q}}{6e\alpha \log(n)} \le \frac{|Q|}{|H|} <
\frac{5}{3\alpha \log(n)}.
\end{displaymath}
\end{enumerate}
\end{theorem}

Note  that  if  $k$  is  odd  and  $k>1$,  then  every  self-conjugate
polynomial  has even  degree.  Hence  when $H$  is one  of  the groups
$\Sp_{2n}(q),  \SO^\pm_{2n}(q)$ or  $\SO_{2n+1}(q)$,  the elements  in
$Q_k$ have  a characteristic  polynomial which has  one self-conjugate
irreducible factor of degree $2k$, and no other irreducible factors of
degree  divisible  by  $k$.   This   follows  from  the  fact  that  a
self-conjugate factor  of degree  divisible by an  odd $k$  would have
degree divisible  by $2k.$ Thus, when  $k$ is odd, for  all the groups
$H$, the  conditions in Definition~\ref{def:Q} on $c_g(x)$  for $g$ to
lie in $Q_k$ are that there is an irreducible factor of degree $\alpha
k$ and  all other  irreducible factors have  degrees not  divisible by
$k$.
For $n\geq 5$  and odd  integers $k  > \log(n)$  we  have in
particular  that  $k\ge  3$  and  therefore,  by  the Theorem of 
\cite{Zsigmondy},
$q^k-1$ has  a primitive  prime divisor.   Moreover, for $n\ge 5$ there
is at least one odd integer $k$ such that $\log(n) < k \le
2\log(n)$. For such $k$,
Theorem~\ref{the:main} yields a  lower bound for  the proportion $|Q_k|/|H|,$
namely   $|Q_k^{\ppd}|/|H|\geq (1-\frac{1}{\alpha\log(n)})\frac{1}{6e\alpha
  \log(n)}  >  \frac{1}{16 e\alpha\log(n)}$,  since  
$(1-\frac{1}{\alpha\log(n)})/6 > 
(1-\frac{1}{\log(5)})/6 >  1/16.$
 Thus we have the following corollary.

\begin{corollary}\label{cor:main2}
Let $q$ be a prime power  and $\F_q$ a finite field with $q$ elements.
Let $H$, $d$, $\alpha$, $\delta$, $\epsilon$ be as in one of the lines
of Table~$\ref{tab:G}$.  If $n \ge 5$ then for each odd integer
$k$ such that $\log(n) < k \le 2\log(n)$, the proportion of elements
in $H$ which have a characteristic
polynomial with  an irreducible  factor of degree  $\alpha k$  and all
other  irreducible  factors of  degree  not divisible by   $k$ is  at  least
$\frac{1}{16 e\alpha\log(n)}$.
\end{corollary}

Moreover, if  $n\geq5$ then, by \cite[Theorem~8.7]{NZM}, there is an odd prime
$k$ satisfying $\log(n) < k < 2 \log(n)$.   For odd primes $k$,
the divisibility condition `not divisible by $k$' is equivalent to the
condition `coprime to $k$'. Thus an immediate consequence of
Corollary~\ref{cor:main2} is the following result.

\begin{corollary}\label{cor:main3}
Let $q$ be a prime power and $\F_q$ a finite field with $q$
Let $H$, $d$, $\alpha$, $\delta$, $\epsilon$ be as in one of the lines
of Table~$\ref{tab:G}$.  If $n \ge 5$ then 
with probability at least $1/(16e \alpha \log(n)), $
the characteristic polynomial of a random element of $H$ has some
irreducible factor, say of degree  
$k$,  such that $\log(n) < k \leq 2 \log(n)$ and the degrees of all
other irreducible factors are coprime to $k$. 
\end{corollary}

\section{The ingredients of the proofs}\label{sec:ingr}

We  quote a version  of a  theorem, originally  employed in  different
contexts independently by \cite{IKS} and \cite{Lehrer}, which
allows us to  estimate the 
proportion of elements  in any of the groups  $H$ of Table~\ref{tab:G}
which lie in one of the sets $Q_k$, $Q_k^{\ppd}$ and $Q_k^{\full}$.  
We use the notation introduced in Section~\ref{sec:ew} for $F, W$.

\begin{theorem}\label{thm:Q}
Let $H$ be one of  the  groups in Table~$\ref{tab:G}$, let $k$ be a 
positive integer with $k\geq\log(n)$, and
let $Q$ denote one of the sets $Q_k$, $Q_k^{\ppd}$ or $Q_k^{\full}.$  Then
\begin{displaymath}
\frac{|Q|}{|H|} = \sum_{C} \frac{|C|}{|W|}\frac{|Q\cap  T_C|}{|T_C|},
\end{displaymath}
where the sum is over all $F$-conjugacy classes $C$ of $W$ and $T_C$
is a  representative of the $H$-conjugacy class of $F$-stable maximal tori
corresponding to the $F$-conjugacy class $C$.
\end{theorem}

\noindent \emph{Proof:}\quad
Clearly $Q$ is closed under conjugation by elements in $G$. Moreover,
for an element $g\in H$ having
Jordan decomposition $g = us$ with $u$ unipotent and $s$ semisimple, we
have $g\in Q$ if and only if $s\in Q,$ as $g$ and $s$ have the
same characteristic polynomial in their actions on the underlying vector space. 
The result follows as in  the proof of
 \cite[Theorem 6.2]{IKS} or the proof of  \cite[Lemma~2.3]{LNP} or from
 \cite[Theorem~1.3]{NP}.

\subsection{Relevant conjugacy classes}
The proof of the main theorem will amount to identifying relevant conjugacy
classes $C$ in $W$ and estimating the proportion of elements of $Q$
that lie in the corresponding torus. The following lemma allows us to 
identify the relevant conjugacy classes when $Q= Q_k$ for some
positive integer $k$. It follows directly from the discussion in 
Section~\ref{sec:ew}.

\begin{lemma}\label{lem:cycles}
Let $H, n$ be as in one of the lines of Table~$\ref{tab:G}$, and let
$2 \le k < n$. Then a maximal $F$-stable torus $T_C$ satisfies $T_C
\cap Q_k \not=\varnothing$ if and only if the corresponding
$F$-conjugacy class $C$ contains elements $w$ for which 
\begin{itemize}
\item[$(a)$] if $H=\SL_{n}(q)$, then $w$ contains a $k$-cycle and all
  other cycle lengths are not divisible by $k$; 
\item[$(b)$] if $H=\SU_{n}(q)$, then $k$ is odd,  $w$ contains a
  $k$-cycle, and all   other cycle lengths are not divisible by $k$;
\item[$(c)$] if
$H=\Sp_{2n}(q),$ $\SO_{2n+1}(q)$, or $\SO^\pm_{2n}(q)$, then  $w$
contains a  unique negative   cycle of length $k$ (that is, its image in $S_n$ has 
length $k$). Moreover, $k$ does not divide the length
of any other cycle, positive or negative.
\end{itemize}
\end{lemma}

\subsection{Proportions of elements in tori which lie in $Q$}

Our aim in this subsection is to find good upper and lower bounds
for $|T_C\cap Q_k|/|T_C|$ for a given $F$-conjugacy class $C$ of $W$.

The following  $O({1/\log(m)})$ lower bound for  $\varphi(m)/m$ is
well known (see for example \cite[Lemma~10.6]{NP98}),
where $\varphi(m)$ is the Euler phi-function, that is, the number 
of positive integers less than $m$ and coprime to $m$.

\begin{lemma}\label{lem:varphi}
For all integers $m$ with $m\ge 3$ there is a constant $a$ such
that  $\varphi(m)/m \ge a \log\log(m)/\log(m)$ and in particular
  $\varphi(m)/m \ge \frac{\log(2)}{\log(m)}$. 
\end{lemma}

\begin{lemma}\label{lem:qk}
Let $C$ be an $F$-conjugacy class of $W$ and $T_C$ a representative of the
corresponding class of $F$-stable maximal tori. Let $k$ be an integer
such that $k\ge \log(n)$ and $9 < q^k$. Then  
\begin{enumerate}
   \item[(a)]   If $T_C\cap Q_k\ne
\emptyset$, then 
\begin{math}1-\frac{2}{q^{k/2}} \le \frac{|T_C\cap
Q_k|}{|T_C| } \le 1.
\end{math}
\item[(b)]  If $T_C\cap
Q^{\ppd}_k\ne \emptyset$, then  
\begin{math}
1-\frac{1}{\alpha k} \le \frac{|T_C\cap Q^{\ppd}_k|}{|T_C| } \le 1  
\end{math} with $\alpha$ as in Table~$\ref{tab:G}$.
\item[(c)] If $T_C\cap
Q^{\full}_k\ne \emptyset$, then  
\begin{math}
\frac{\log(2) }{k\log(q+1)} \le \frac{|T_C\cap
Q_k^{\full}|}{|T_C| } \le 1.   
\end{math}
\end{enumerate}

\end{lemma}

\noindent \emph{Proof:}\quad
Clearly the upper bound holds in all three cases.
Statement (a) is proved in  the proof of \cite[Lemma~4.3(a)]{NPP}.

To prove part (b) let $m =\alpha k$,  with $\alpha$ as in Table~$\ref{tab:G}$. 
Since $T_C\cap Q^{\ppd}_k\ne \emptyset$ the abelian group 
 $T_C$ contains elements of order divisible by a primitive prime
divisor $r$ of $q^{\delta m}-1$. All elements of $T_C$ lying outside a subgroup of index $r$ have order divisible by $r$, and hence
the proportion of such elements in $T_C$ is at least $1-\frac{1}{r}
> 1-\frac{1}{m}$, as $r\geq m+1$. Thus the result follows.

For part (c), note that the $F$-stable maximal torus
$T_C$ contains a direct factor $Z$ of order
$\ell:=q^k-1$ or $q^k+1$. If $t\in T_C$ is such that its projection 
$\pi(t)$ onto $Z$ generates $Z$, then $t$ lies in $Q^{\full}_k$ since
the $r$-part $r^a$ of $q^k-1$ or $q^{2k}-1$, respectively, divides $\ell$,
and hence divides $|\pi(t)|$.
The number of generators of $\pi(T_C)$ is at least
$\varphi(q^k-1)$ or $\varphi(q^k+1)$, respectively. Hence a lower bound for
$|T_C\cap Q_k^\full|/|T_C|$ is
$\varphi(\ell)/\ell$.
By Lemma~\ref{lem:varphi}, $\varphi(\ell)/\ell\ge 
\log(2)/\log(\ell)$ which is at least $\log(2)/k\log(q+1)$ since 
$\log(\ell)\leq \log(q^k+1)\leq \log((q+1)^k)$.

\subsection{Proportions of elements in symmetric groups}

Let  $n, m$  be positive  integers with  $m <  n$, and  let $b_{m}(n)$
denote the  proportion of elements  in the symmetric group  $S_n$ with
exactly  one cycle  of  length $m$  and  all other  cycle lengths  not
divisible  by $m$. 

\begin{lemma}\label{lem:sn}
If $m\ge3$ and $\log(n) \le  m \le n-m$, then
\begin{math}
\frac{1}{3e m}  \le b_{m}(n) \le \frac{5}{3m}.
\end{math}

\end{lemma}

\noindent \emph{Proof:}\quad
Let $p_{\neg m}(n-m)$ denote the proportion of elements in $S_{n-m}$ which
contain  no cycle of length divisible by $m$. 
Observe first that 
$b_{m}(n)  = \frac{1}{m} p_{\neg   m}(n-m).$ Putting 
$c(j) = \frac{1}{\Gamma(1-1/j)}$  and applying the inequality
derived from  \cite[Theorem 2.3(b)] {Bealsetal} in 
the proof of  \cite[Lemma 4.2]{LNP}  we
have
$$\frac{m^{1/m}c(m)}{m(n- m)^{1/m}}
\left(1-\frac{1}{n-m}\right) \le b_{m}(n) \le \frac{m^{1/m} c(m) }
{m(n- m)^{1/m}}\left(1 + \frac{2}{n- m}\right) .$$ 

By our hypothesis $n - m\ge 3$, and also $1/2 \le c(m) \le 1$.
Hence
$$\frac{1}{3} \frac{m^{1/m}}{m(n-m)^{1/m}}
\le b_{m}(n) \le \frac{5}{3}\frac{m^{1/m}}{m(n-  m)^{1/m}}  .$$ 

Now \begin{math}
      \left(\frac{m}{n-m}\right)^{1/m}\leq 1  
    \end{math}
 since $m\leq n-m$. This gives the required upper bound. We claim that
 \begin{math}
      \left(\frac{m}{n-m}\right)^{1/m}\geq \frac{1}{e},  
    \end{math}
or equivalently, $e^m\geq (n-m)/m$. To see that this is true, observe 
that, since $m\geq \log(n)$, we have $e^m\geq e^{\log(n)} =n>(n-m)/m$.
This proves the claim and yields the required lower bound.

\subsection{Proof of Theorem~\ref{the:main}}\label{sub:main}

Let $Q$ denote one of the sets  $Q_k$, $Q_k^{\ppd}$ or $Q_k^{\full}$,
where $k\ge\log(n)$.
By Theorem~\ref{thm:Q},
$$
\frac{|Q|}{|H|} = \sum_{C} \frac{|C|}{|W|}\frac{|Q\cap  T_C|}{|T_C|},
$$
where the sum is over all $F$-conjugacy classes $C$ of $W$, and
$T_C$ is a  representative of the $H$-conjugacy class of $F$-stable
maximal tori 
corresponding to $C$ such that 
$T_C\cap Q\ne \emptyset$. 
By 
Lemma~\ref{lem:qk}, we see that $\ell_{k,Q}$, as defined in Table~\ref{tab:ell},
is a uniform lower bound for each non-zero $|T_C\cap Q|/|T_C|$ and we
use $1$ as a uniform upper bound for this quantity.

In  lines 1 and 2 of  Table~\ref{tab:G}, the proportion
$\sum_{C} |C|/|W|$ is equal to $b_k(n)$, as defined in Section~4.3. 
In  lines 3--5 of Table~\ref{tab:G}, 
the $F$-conjugacy classes $C$ which correspond to tori containing 
elements of $Q_k$ are those with
one negative cycle of length $k$ and  no other 
cycles of length divisible by $k$. 
By \cite[Lemma 4.2(c,d)]{LNP} it
follows that $\sum_{C} |C|/|W| =  b_{k}(n)/2.$ Thus in all cases
 $\sum_{C} |C|/|W| =  b_{k}(n)/\alpha,$ with $\alpha$ as in Table~\ref{tab:G}.
The upper and lower bounds for part $(a)$ now follow, since by 
Lemma~\ref{lem:sn} $1/(3e\alpha k) \leq b_n(k)< 5/(3\alpha k)$.

 To see that part (b) holds, 
 assume
 that $\log(n)< k \leq 2\log(n)$. The upper bound follows immediately
 from (a) since $k>\log(n)$. Similarly the lower 
 bound follows using the bounds on $k$  to obtain the lower bound
 $m_{n,Q}$ for $\ell_{k,Q}$ with $k$ in the given range.

\subsection{Proof of Theorem~\ref{the:short-main}}

The upper  bound follows immediately  from Theorem~\ref{the:main}. For
the lower bound for $Q_k$  and $Q_k^\ppd$ we may choose $Q= Q_k^\ppd$.
Observe that  the quantity  $\ell_{k,Q}$ in Table~\ref{tab:ell}  is at
least $\ell_{k,Q}  \ge 1-1/(\alpha k)\geq 2/3$  since $\alpha k\geq3$.
Then by Theorem~\ref{the:main}, $|Q|/|H|\geq 2/(9e\alpha k)$. The
lower bound for $Q_k^\full$  follows immediately from 
Theorem~\ref{the:main}.

\section*{Acknowledgements}

The authors express their sadness at the recent passing of their
colleague \'Akos Seress, and wish to acknow\-ledge his pivotal role in
formulating the problem addressed in this paper, and indeed  
their appreciation of working with him over a number of years. 
The paper was originally submitted in December 2011 for the issue
celebrating the sixtieth birthday of L\'aszl\'o Babai, and the authors
extend very belated congratulations to Laci.  The authors also thank
Frank L\"ubeck for many  valuable suggestions on
the  exposition.   The  research  forms part  of  Australian  Research
Council Discovery Grants DP110101153 and DP140100416.   The first
author acknowledges a DFG grant  in SPP1489.   During the writing of
this paper the second author  was supported  by Australian Research
Council Federation Fellowship FF0776186. 

\nocite{*}
\bibliographystyle{abbrvnat}
\def\cprime{$'$}

{\tiny
\begin{enumerate}
\item[1] Centre for the Mathematics of Symmetry and Computation,\\
School  of Mathematics and Statistics M019,\\
The University of Western Australia,\\
35 Stirling Highway,\\
Crawley, WA 6009,\\
Australia
\item[2]
Lehrstuhl D f\"ur Mathematik,\\
RWTH Aachen University,\\
Templergraben 64, \\
52064 Aachen,\\
Germany
\item[3]
King Abdulazziz University, \\
Jeddah,\\
Saudi Arabia\\
\end{enumerate}
{\tt Cheryl.Praeger@uwa.edu.au, Alice.Niemeyer@uwa.edu.au}
}
\end{document}